\documentclass[12pt]{article}%
\usepackage{amsmath}
\usepackage{amsfonts}
\usepackage{amssymb}
\usepackage{graphicx}%
\providecommand{\U}[1]{\protect\rule{.1in}{.1in}}
\newtheorem{theorem}{Theorem}
\newtheorem{acknowledgement}[theorem]{Acknowledgement}

\newtheorem{corollary}[theorem]{Corollary}

\newtheorem{example}[theorem]{Example}

\newenvironment{proof}[1][Proof]{\noindent\textbf{#1.} }{\ \rule{0.5em}{0.5em}}

\begin{document}

\title{Polynomial solutions of nonlinear integral equations}
\author{Diego Dominici\\Department of Mathematics\\State University of New York at New Paltz\\1 Hawk Dr. Suite 9\\New Paltz, NY 12561-2443\\dominicd@newpaltz.edu}
\maketitle

\begin{abstract}
We analyze the polynomial solutions of a nonlinear integral equation, generalizing the work of C. Bender and E. Ben-Naim. We show that, in some cases, an orthogonal solution exists and we give its general form in terms of kernel polynomials.

\end{abstract}



\section{Introduction}

In \cite{MR2304926}, C. Bender and E. Ben-Naim studied the polynomial
solutions of the nonlinear integral equation%
\begin{equation}
\int\limits_{a}^{b}P(y)P\left(  x+y\right)  \omega(y)dy=P(x). \label{eq1}%
\end{equation}
They showed that the solutions $P_{n}(x)$ are orthogonal with respect to the
measure $y\omega(y)$ and considered other equations of the form%

\begin{equation}
\int\limits_{a}^{b}P(y)P\left[  F(x,y)\right]  \omega(y)dy=P(x), \label{eq0}%
\end{equation}
with
\[
F(x,y)=xy,\quad x+a_{1}+a_{2}y\text{ \ \ and \ \ }x+f(y).
\]

The purpose of this paper is to generalize their results to the case
$F(x,y)=\alpha(y)+x\beta(y),$ for arbitrary functions $\alpha(y)$ and
$\beta(y).$ We also try to understand the nature of the families of orthogonal
polynomials that arise as solutions of (\ref{eq1}) and show that they are in
fact kernel polynomials associated with the weight $\omega(y).$

\section{General case}

Let $\omega(y)$ be a non-negative integrable function on the interval $\left(
a,b\right)  ,$ such that
\begin{equation}
\int\limits_{a}^{b}\omega(y)dy=1\label{omega}%
\end{equation}
and let $%
\mathcal{L}%
_{\omega}$ be the linear functional defined by
\begin{equation}%
\mathcal{L}%
_{\omega}\left[  f\right]  =\int\limits_{a}^{b}f\left(  y\right)
\omega(y)dy.\label{L}%
\end{equation}
We say that a sequence of polynomials $\left(  P_{n}\right)  $ is an
orthogonal polynomial sequence (OPS) with respect to $%
\mathcal{L}%
_{\omega}$ if \cite{MR0481884}:

\begin{enumerate}
\item $P_{n}(x)$ is a polynomial of degree $n.$

\item $%
\mathcal{L}%
_{\omega}\left[  P_{n}P_{m}\right]  =h_{n}\delta_{n,m},\quad n,m=0,1,\ldots,$

where $h_{n}\neq0$ for all $n$ and $\delta_{n,m}$ is Kronecker's delta.
\end{enumerate}

To warranty the existence of a polynomial sequence solution $\left(
P_{n}\right)  $, we consider the special form of equation (\ref{eq0})%

\begin{equation}
\int\limits_{a}^{b}P_{n}(y)P_{n}\left[  \alpha(y)+x\beta(y)\right]
\omega(y)dy=P_{n}(x), \label{eq}%
\end{equation}
where $\alpha(y)$ and $\beta(y)$ are integrable functions on $\left(
a,b\right)  .$

\begin{example}
Let
\[
\omega(y)=\frac{3}{2}y^{2},\quad a=-1,b=1,\quad\alpha(y)=\frac{5}{3}%
y,\quad\beta(y)=\mu\neq0.
\]
Then, we have%
\[
P_{0}(x)=1,\quad P_{1}(x)=\frac{1}{\mu}\pm\frac{\sqrt{\mu-1}}{\mu}%
x,\quad\ldots.
\]

\end{example}

\begin{example}
Let
\[
\omega(y)=\frac{3}{2}y^{2},\quad a=-1,b=1,\quad\alpha(y)=\frac{3}{20}\mu
y,\quad\beta(y)=y.
\]
Then, we have%
\[
P_{0}(x)=1,\quad P_{1}(x)=\frac{1\pm\sqrt{1-\mu}}{2}+\frac{5}{3}x,\quad
\ldots.
\]

\end{example}

The previous examples illustrate how, even for simple functions, the integral
equation (\ref{eq}) can have unique or multiple solutions which are real or
complex depending on the choice of the parameter $\mu$.

Writing
\begin{equation}
P_{n}(x)=\sum\limits_{k=0}^{n}c_{k}x^{k}, \label{Pn}%
\end{equation}
we have%
\begin{equation}
P_{n}\left[  \alpha(y)+x\beta(y)\right]  =\sum\limits_{k=0}^{n}\gamma
_{k}\left(  y\right)  x^{k}, \label{Pn1}%
\end{equation}
where%
\[
\gamma_{k}\left(  y\right)  =\beta^{k}(y)\sum\limits_{j=k}^{n}c_{j}\binom
{j}{k}\alpha^{j-k}(y).
\]
Using (\ref{Pn}) and (\ref{Pn1}) in (\ref{eq}), we get%
\begin{equation}%
\mathcal{L}%
_{\omega}\left[  P_{n}\gamma_{k}\right]  =c_{k},\quad0\leq k\leq n.
\label{matrix1}%
\end{equation}

Introducing the matrix $\mathbf{A}$ defined by%
\begin{equation}
\mathbf{A}=%
\begin{bmatrix}
\binom{0}{0}%
\mathcal{L}%
_{\omega}\left[  P_{n}\right]  & \binom{1}{0}%
\mathcal{L}%
_{\omega}\left[  P_{n}\alpha\right]  & \binom{2}{0}%
\mathcal{L}%
_{\omega}\left[  P_{n}\alpha^{2}\right]  & \cdots & \binom{n}{0}%
\mathcal{L}%
_{\omega}\left[  P_{n}\alpha^{n}\right] \\
0 & \binom{1}{1}%
\mathcal{L}%
_{\omega}\left[  P_{n}\beta\right]  & \binom{2}{1}%
\mathcal{L}%
_{\omega}\left[  P_{n}\alpha\beta\right]  & \cdots & \binom{n}{1}%
\mathcal{L}%
_{\omega}\left[  P_{n}\alpha^{n-1}\beta\right] \\
0 & 0 & \binom{2}{2}%
\mathcal{L}%
_{\omega}\left[  P_{n}\beta^{2}\right]  & \cdots & \binom{n}{2}%
\mathcal{L}%
_{\omega}\left[  P_{n}\alpha^{n-2}\beta^{2}\right] \\
\vdots & \vdots & 0 & \ddots & \vdots\\
0 & 0 & 0 & \cdots & \binom{n}{n}%
\mathcal{L}%
_{\omega}\left[  P_{n}\beta^{n}\right]
\end{bmatrix}
\label{A}%
\end{equation}
and the vector
\[
\mathbf{C}^{T}=%
\begin{bmatrix}
c_{0}, & \cdots & c_{n}%
\end{bmatrix}
,
\]
we see from (\ref{matrix1}) that $\mathbf{C}$ is an eigenvector of
$\mathbf{A}$ with corresponding eigenvalue $1.$

Therefore, to have a solution $\mathbf{C}$ different from the zero vector, it
must be true that $%
\mathcal{L}%
_{\omega}\left[  P_{n}\beta^{k}\right]  =1$ for some $0\leq k\leq n.$ Note
that if we impose the condition $%
\mathcal{L}%
_{\omega}\left[  P_{n}\beta^{n}\right]  =1,$ then the vector $\mathbf{C}%
^{T}\mathbf{=}%
\begin{bmatrix}
0, & \cdots & 0, & c_{n}%
\end{bmatrix}
$ is always an eigenvector of $\mathbf{A}$. However, this leads to trivial
sequences of the form $P_{n}(x)=c_{n}x^{n}.$

A possible non-trivial solution of the equation $\mathbf{AC}=\mathbf{C}$ is to
take $\mathbf{A=I,}$ where $\mathbf{I}$ denotes the identity matrix. Thus, we
require that%
\begin{equation}%
\mathcal{L}%
_{\omega}\left[  P_{n}\alpha^{j-i}\beta^{i}\right]  =\delta_{i,j},\quad0\leq
i\leq j\leq n. \label{sys}%
\end{equation}
Since (\ref{sys}) is a system of $\binom{n+2}{2}$ equations with $n+1$
unknowns, it admits (if any) infinitely many solutions. In order to have a
unique solution, we consider the following cases:

\begin{enumerate}
\item $\alpha(y)=0,\quad\beta(y)\neq1.$

We see from (\ref{A}) that for $\mathbf{A}$ to be equal to the identity
matrix, we need to have%
\begin{equation}%
\mathcal{L}%
_{\omega}\left[  P_{n}\beta^{i}\right]  =1,\quad0\leq i\leq n, \label{cond1}%
\end{equation}
which, using (\ref{omega}), we can rewrite as%
\[%
\mathcal{L}%
_{\omega}\left[  P_{n}\left(  1-\beta^{i}\right)  \right]  =0,\quad0\leq i\leq
n,
\]
or%
\[%
\mathcal{L}%
_{\left(  \beta-1\right)  \omega}\left[  P_{n}\frac{\left(  \beta
^{i}-1\right)  }{\beta-1}\right]  =0,\quad0\leq i\leq n.
\]
Thus, (\ref{cond1}) is equivalent to%
\begin{equation}%
\mathcal{L}%
_{\omega}\left[  P_{n}\right]  =1,\quad%
\mathcal{L}%
_{\left(  \beta-1\right)  \omega}\left[  P_{n}\beta^{i}\right]  =0,\quad0\leq
i\leq n-1. \label{orthbeta}%
\end{equation}

If $\beta(y)$ is linear, it follows from (\ref{orthbeta}) that $\left(
P_{n}\right)  $ will be a sequence of orthogonal polynomials with respect to
the linear functional $%
\mathcal{L}%
_{\left(  \beta-1\right)  \omega},$ provided that $%
\mathcal{L}%
_{\left(  \beta-1\right)  \omega}\left[  P_{n}\beta^{n}\right]  \neq0.$ For
this last condition to be true, $\beta\left(  y\right)  -1$ must not vanish in
the interval $\left(  a,b\right)  .$ Hence, $\beta\left(  y\right)  $ should
be of the form%
\[
\beta\left(  y\right)  =\sigma\left(  y-\zeta\right)  +1,
\]
with $\sigma\neq0$ and $\zeta\notin\left(  a,b\right)  .$

\item $\alpha(y)\neq0,\quad\beta(y)=1.$

In this case, we must impose that%
\[%
\mathcal{L}%
_{\omega}\left[  P_{n}\right]  =1,\quad%
\mathcal{L}%
_{\omega}\left[  P_{n}\alpha^{i}\right]  =0,\quad1\leq i\leq n,
\]
or, equivalently,%
\begin{equation}%
\mathcal{L}%
_{\omega}\left[  P_{n}\right]  =1,\quad%
\mathcal{L}%
_{\alpha\omega}\left[  P_{n}\alpha^{i}\right]  =0,\quad0\leq i\leq n-1.
\label{orthalpha}%
\end{equation}
If
\[
\alpha(y)=\tau\left(  y-\varsigma\right)  ,
\]
with $\tau\neq0$ and $\varsigma\notin\left(  a,b\right)  ,$ the polynomials
$P_{n}(x)$ will be orthogonal with respect to $%
\mathcal{L}%
_{\alpha\omega}.$
\end{enumerate}

We summarize the results of this section in the following theorems.

\begin{theorem}
Let $\beta(y)\neq1$ on $\left(  a,b\right)  $ and suppose that $\Delta_{n}%
\neq0$ for all $n,$ with%
\begin{equation}
\Delta_{n}=%
\begin{vmatrix}
1 &
\mathcal{L}%
_{\omega}\left[  y\right]  & \cdots &
\mathcal{L}%
_{\omega}\left[  y^{n}\right] \\%
\mathcal{L}%
_{\left(  \beta-1\right)  \omega}\left[  1\right]  &
\mathcal{L}%
_{\left(  \beta-1\right)  \omega}\left[  y\right]  & \cdots &
\mathcal{L}%
_{\left(  \beta-1\right)  \omega}\left[  y^{n}\right] \\
\vdots & \vdots & \ddots & \vdots\\%
\mathcal{L}%
_{\left(  \beta-1\right)  \omega}\left[  \beta^{n-1}\right]  &
\mathcal{L}%
_{\left(  \beta-1\right)  \omega}\left[  y\beta^{n-1}\right]  & \cdots &
\mathcal{L}%
_{\left(  \beta-1\right)  \omega}\left[  y^{n}\beta^{n-1}\right]
\end{vmatrix}
. \label{det}%
\end{equation}
If $\left(  P_{n}\right)  $ is defined by
\begin{equation}
P_{n}(x)=\frac{1}{\Delta_{n}}%
\begin{vmatrix}
1 & x & \cdots & x^{n}\\%
\mathcal{L}%
_{\left(  \beta-1\right)  \omega}\left[  1\right]  &
\mathcal{L}%
_{\left(  \beta-1\right)  \omega}\left[  y\right]  & \cdots &
\mathcal{L}%
_{\left(  \beta-1\right)  \omega}\left[  y^{n}\right] \\
\vdots & \vdots & \ddots & \vdots\\%
\mathcal{L}%
_{\left(  \beta-1\right)  \omega}\left[  \beta^{n-1}\right]  &
\mathcal{L}%
_{\left(  \beta-1\right)  \omega}\left[  y\beta^{n-1}\right]  & \cdots &
\mathcal{L}%
_{\left(  \beta-1\right)  \omega}\left[  y^{n}\beta^{n-1}\right]
\end{vmatrix}
, \label{Pn2}%
\end{equation}
then,%
\[
\int\limits_{a}^{b}P_{n}(y)P_{n}\left[  \beta(y)x\right]  \omega
(y)dy=P_{n}(x)
\]
for all $n.$
\end{theorem}

\begin{proof}
It is clear from (\ref{det}) and (\ref{Pn2}) that
\[%
\mathcal{L}%
_{\omega}\left[  P_{n}\right]  =1,\quad%
\mathcal{L}%
_{\left(  \beta-1\right)  \omega}\left[  P_{n}\beta^{i}\right]  =0,\quad0\leq
i\leq n-1
\]
for all $n\geq1.$ We have%
\begin{align*}
&  \int\limits_{a}^{b}P_{n}(y)P_{n}\left[  \beta(y)x\right]  \omega
(y)dy-P_{n}(x)\\
&  =\int\limits_{a}^{b}P_{n}(y)P_{n}\left[  \beta(y)x\right]  \omega
(y)dy-\int\limits_{a}^{b}P_{n}(y)P_{n}(x)\omega(y)dy\\
&  =\int\limits_{a}^{b}P_{n}(y)\left\{  P_{n}\left[  \beta(y)x\right]
-P_{n}(x)\right\}  \omega(y)dy.
\end{align*}
Using (\ref{Pn}), we get%
\begin{align*}
&  \int\limits_{a}^{b}P_{n}(y)P_{n}\left[  \beta(y)x\right]  \omega
(y)dy-P_{n}(x)\\
&  =\sum\limits_{k=1}^{n}c_{k}%
\mathcal{L}%
_{\omega}\left[  P_{n}\left(  \beta^{k}-1\right)  \right]  x^{k}%
=\sum\limits_{k=1}^{n}c_{k}%
\mathcal{L}%
_{\left(  \beta-1\right)  \omega}\left[  P_{n}\frac{\left(  \beta
^{k}-1\right)  }{\beta-1}\right]  x^{k}\\
&  =\sum\limits_{k=1}^{n}c_{k}\left[  \sum\limits_{j=0}^{k-1}%
\mathcal{L}%
_{\left(  \beta-1\right)  \omega}\left[  P_{n}\beta^{j}\right]  \right]
x^{k}=0.
\end{align*}

\end{proof}

\begin{theorem}
Let $\alpha(y)\neq0$ on $\left(  a,b\right)  $ and suppose that $\Delta
_{n}\neq0$ for all $n,$ with%
\begin{equation}
\Delta_{n}=%
\begin{vmatrix}
1 &
\mathcal{L}%
_{\omega}\left[  y\right]  & \cdots &
\mathcal{L}%
_{\omega}\left[  y^{n}\right] \\%
\mathcal{L}%
_{\alpha\omega}\left[  1\right]  &
\mathcal{L}%
_{\alpha\omega}\left[  y\right]  & \cdots &
\mathcal{L}%
_{\alpha\omega}\left[  y^{n}\right] \\
\vdots & \vdots & \ddots & \vdots\\%
\mathcal{L}%
_{\alpha\omega}\left[  \alpha^{n-1}\right]  &
\mathcal{L}%
_{\alpha\omega}\left[  y\alpha^{n-1}\right]  & \cdots &
\mathcal{L}%
_{\alpha\omega}\left[  y^{n}\alpha^{n-1}\right]
\end{vmatrix}
. \label{det1}%
\end{equation}
If $\left(  P_{n}\right)  $ is defined by
\begin{equation}
P_{n}(x)=\frac{1}{\Delta_{n}}%
\begin{vmatrix}
1 & x & \cdots & x^{n}\\%
\mathcal{L}%
_{\alpha\omega}\left[  1\right]  &
\mathcal{L}%
_{\alpha\omega}\left[  y\right]  & \cdots &
\mathcal{L}%
_{\alpha\omega}\left[  y^{n}\right] \\
\vdots & \vdots & \ddots & \vdots\\%
\mathcal{L}%
_{\alpha\omega}\left[  \alpha^{n-1}\right]  &
\mathcal{L}%
_{\alpha\omega}\left[  y\alpha^{n-1}\right]  & \cdots &
\mathcal{L}%
_{\alpha\omega}\left[  y^{n}\alpha^{n-1}\right]
\end{vmatrix}
, \label{Pn3}%
\end{equation}
then,%
\[
\int\limits_{a}^{b}P_{n}(y)P_{n}\left[  \alpha(y)+x\right]  \omega
(y)dy=P_{n}(x)
\]
for all $n.$
\end{theorem}

\begin{proof}
It is clear from (\ref{det}) and (\ref{Pn2}) that
\[%
\mathcal{L}%
_{\omega}\left[  P_{n}\right]  =1,\quad%
\mathcal{L}%
_{\alpha\omega}\left[  P_{n}\alpha^{i}\right]  =0,\quad0\leq i\leq n-1
\]
for all $n\geq1.$ We have%
\begin{align*}
&  \int\limits_{a}^{b}P_{n}(y)P_{n}\left[  \alpha(y)+x\right]  \omega
(y)dy-P_{n}(x)\\
&  =\int\limits_{a}^{b}P_{n}(y)\left\{  P_{n}\left[  \alpha(y)+x\right]
-P_{n}(x)\right\}  \omega(y)dy\\
&  =\sum\limits_{k=1}^{n}q_{k}(x)%
\mathcal{L}%
_{\omega}\left[  P_{n}\alpha^{k}\right]  =\sum\limits_{k=0}^{n-1}q_{k+1}(x)%
\mathcal{L}%
_{\alpha\omega}\left[  P_{n}\alpha^{k}\right]  =0,
\end{align*}
where we have used (\ref{Pn}) and the polynomials $q_{k}(x)$ are defined by
\[
q_{k}(x)=\sum\limits_{j=k}^{n}c_{j}\binom{j}{k}x^{j-k}.
\]

\end{proof}

\begin{theorem}
\label{Th1}Let $\zeta\notin(a,b)$ and $\left(  P_{n}\right)  $ be an OPS for $%
\mathcal{L}%
_{\left(  y-\zeta\right)  \omega}$ satisfying
\begin{equation}%
\mathcal{L}%
_{\omega}\left[  P_{n}\right]  =1,\quad n=0,1,\ldots. \label{norm}%
\end{equation}
Then,
\begin{equation}
\int\limits_{a}^{b}P_{n}(y)P_{n}\left[  \left(  y-\zeta\right)  \left(
\tau+\sigma x\right)  +x\right]  \omega(y)dy=P_{n}(x). \label{eq3}%
\end{equation}

\end{theorem}

\begin{proof}
Using (\ref{norm}), we see that%
\begin{align*}
&  \int\limits_{a}^{b}P_{n}(y)P_{n}\left[  \left(  y-\zeta\right)  \left(
\tau+\sigma x\right)  +x\right]  \omega(y)dy-P_{n}(x)\\
&  =\int\limits_{a}^{b}P_{n}(y)\left\{  P_{n}\left[  \left(  y-\zeta\right)
\left(  \tau+\sigma x\right)  +x\right]  -P_{n}(x)\right\}  \omega(y)dy\\
&  =\sum\limits_{k=1}^{n}q_{k}(x)\left(  \tau+\sigma x\right)  ^{k}%
\mathcal{L}%
_{\omega}\left[  \left(  y-\zeta\right)  ^{k}P_{n}\right] \\
&  =\sum\limits_{k=0}^{n-1}q_{k+1}(x)\left(  \tau+\sigma x\right)  ^{k+1}%
\mathcal{L}%
_{\left(  y-\zeta\right)  \omega}\left[  \left(  y-\zeta\right)  ^{k}%
P_{n}\right]  =0.
\end{align*}

\end{proof}

The two main cases considered in \cite{MR2304926}, correspond to the special
values
\begin{align*}
\zeta &  =0,\tau=1,\sigma=0,\quad\text{for \ }P_{n}(y)P_{n}\left(  x+y\right)
\\
\zeta &  =1,\tau=0,\sigma=1,\quad\text{for \ }P_{n}(y)P_{n}\left(  xy\right)
.
\end{align*}

Unfortunately, the reciprocal of Theorem \ref{Th1} is not true in general. For
example, taking%
\[
\omega(y)=e^{-y},\quad a=0,\quad b=\infty,\quad\zeta=0,\quad\tau=1,\quad
\sigma=1,
\]
we get as a possible solution of (\ref{eq3})%

\begin{align*}
P_{0}(x)  &  =1,\quad P_{1}(x)=2-x,\quad P_{2}(x)=\frac{7}{5}-\frac{1}%
{5}x-\frac{1}{10}x^{2},\\
P_{3}(x)  &  =\frac{43}{17}-\frac{32}{17}x+\frac{3}{34}x^{2}+\frac{1}{34}%
x^{3},\ldots.
\end{align*}
We have
\[%
\mathcal{L}%
_{y\omega}\left[  P_{1}\right]  =0,\quad%
\mathcal{L}%
_{y\omega}\left[  P_{2}\right]  =\frac{2}{5},\quad%
\mathcal{L}%
_{y\omega}\left[  P_{3}\right]  =0,\quad%
\mathcal{L}%
_{y\omega}\left[  yP_{3}\right]  =-\frac{10}{17}\ldots\quad
\]
and therefore $\left(  P_{n}\right)  $ is not an OPS for $%
\mathcal{L}%
_{\left(  y-\zeta\right)  \omega}.$

In the next section, we shall see that the only solutions of (\ref{eq3}) which
are an OPS for $%
\mathcal{L}%
_{\left(  y-\zeta\right)  \omega},$ consist of the so called kernel
polynomials corresponding to $%
\mathcal{L}%
_{\omega}.$

\section{Kernel polynomials}

Let $\left(  \mathfrak{p}_{n}\right)  $ be the sequence of orthonormal
polynomials with respect to $%
\mathcal{L}%
_{\omega}$ defined by (\ref{L}). The kernel polynomials $K_{n}(x;\zeta)$
\ corresponding to $%
\mathcal{L}%
_{\omega}$ with parameter $\zeta$ are defined by \cite{MR0372517}%
\begin{equation}
K_{n}(x;\zeta)=\sum\limits_{k=0}^{n}\frac{\mathfrak{p}_{k}\left(
\zeta\right)  }{%
\mathcal{L}%
_{\omega}\left[  \mathfrak{p}_{k}^{2}\right]  }\mathfrak{p}_{k}\left(
x\right)  , \label{kernel1}%
\end{equation}
where $\mathfrak{p}_{n}\left(  \zeta\right)  \neq0$ for all $n.$ Using the
Christoffel-Darboux Identity \cite{MR1688958}, we have%
\[
K_{n}(x;\zeta)=\frac{1}{%
\mathcal{L}%
_{\omega}\left[  \mathfrak{p}_{n}^{2}\right]  }\frac{\mathfrak{p}_{n+1}\left(
x\right)  \mathfrak{p}_{n}\left(  \zeta\right)  -\mathfrak{p}_{n}\left(
x\right)  \mathfrak{p}_{n+1}\left(  \zeta\right)  }{x-\zeta}.
\]
The kernel polynomials $K_{n}(x;\zeta)$ have the following properties
\cite{MR0481884}:

\begin{enumerate}
\item They are orthogonal with respect to the functional $%
\mathcal{L}%
_{\left(  y-\zeta\right)  \omega}.$

\item They have the reproducing property%
\[%
\mathcal{L}%
_{\omega}\left[  K_{n}(y;\zeta)p_{n}\left(  y\right)  \right]  =p_{n}\left(
\zeta\right)  ,
\]
for any polynomial $p_{n}\left(  x\right)  $ of degree less or equal than $n.$
\end{enumerate}

It follows that, up to a multiplicative constant $\lambda$, the kernel
polynomials $K_{n}(x;\zeta)$ are solutions of (\ref{eq1}). To find $\lambda$
we use (\ref{norm}) and obtain%
\[
1=%
\mathcal{L}%
_{\omega}\left[  \lambda K_{n}(x;\zeta)\right]  =\lambda\sum\limits_{k=0}%
^{n}\frac{\mathfrak{p}_{k}\left(  \zeta\right)  }{%
\mathcal{L}%
_{\omega}\left[  \mathfrak{p}_{k}^{2}\right]  }%
\mathcal{L}%
_{\omega}\left[  \mathfrak{p}_{k}\left(  x\right)  \right]  =\lambda.
\]
Thus, we have the following result.

\begin{corollary}
The only OPS $\left(  P_{n}\right)  $ which is a solution of the nonlinear
integral equation
\[
\int\limits_{a}^{b}P_{n}(y)P_{n}\left[  \left(  y-\zeta\right)  \left(
\tau+\sigma x\right)  +x\right]  \omega(y)dy=P_{n}(x)
\]
is $P_{n}(x)=K_{n}(x;\zeta),$ where $K_{n}(x;\zeta)$ is defined by
(\ref{kernel1}).
\end{corollary}

\begin{example}
Let
\[
\omega(y)=\frac{1}{2},\quad a=-1,\quad b=1,\quad\zeta=1,\quad\tau
=1,\quad\sigma=0.
\]
Then, we have%
\[
P_{0}(x)=1,\quad P_{1}(x)=1+3x,\quad P_{2}(x)=-\frac{3}{2}+3x+\frac{15}%
{2}x^{2},\quad\ldots.
\]
If we denote by $\mathbf{P}_{n}(x)$ the Legendre polynomials, defined by
\cite{MR2191786}%
\[
\mathbf{P}_{n}(x)=\sum\limits_{k=0}^{n}\binom{n}{k}\binom{-n-1}{k}\left(
\frac{1-x}{2}\right)  ^{k},
\]
then it follows from (\ref{kernel1}) that%
\[
P_{n}(x)=\sum\limits_{k=0}^{n}\frac{\mathbf{P}_{k}(1)}{\left(  2k+1\right)
^{-1}}\mathbf{P}_{k}(x)=\sum\limits_{k=0}^{n}\left(  2k+1\right)
\mathbf{P}_{k}(x).
\]

\end{example}

\begin{example}
Again, let%
\[
\omega(y)=e^{-y},\quad a=0,\quad b=\infty,\quad\zeta=0,\quad\tau=1,\quad
\sigma=1.
\]
Then,
\[
P_{n}(x)=\sum\limits_{k=0}^{n}L_{k}(0)L_{k}(x)=\sum\limits_{k=0}^{n}L_{k}(x),
\]
$\left(  P_{n}\right)  $ is an OPS for $%
\mathcal{L}%
_{\left(  y-\zeta\right)  \omega},$ where%
\[
L_{n}(x)=\sum\limits_{k=0}^{n}\frac{1}{k!}\binom{n}{k}\left(  -x\right)  ^{k}%
\]
denotes the Laguerre polynomial \cite{MR2191786}. We have%
\begin{align*}
P_{0}(x)  &  =1,\quad P_{1}(x)=2-x,\quad P_{2}(x)=3-3x+\frac{1}{2}x^{2},\\
P_{3}(x)  &  =4-6x+2x^{2}-\frac{1}{6}x^{3},\ldots.
\end{align*}

\end{example}

\section{Concluding remarks}

We have studied the polynomial solutions of the nonlinear integral equation
(\ref{eq}). We have shown that, in some cases, a solution which is an OPS
exists and we have given the general form of these orthogonal solutions.

However, much remains to be discovered about the solutions of (\ref{eq}). A
few questions that come to mind are:

\begin{enumerate}
\item For which choice of $\alpha$ and $\beta$ will there be a unique solution?

\item Is it possible to describe all possible solutions?

\item For which values of $\zeta,\tau$ and $\sigma$ will the solution of
(\ref{eq3}) be unique? It seems that for this to be true, one needs to
consider the symmetric case, when
\[
\zeta\sigma+\tau=1.
\]
Is this condition sufficient?
\end{enumerate}

We hope that other researchers will find this problem interesting and continue
its analysis.

\begin{acknowledgement}
\vspace*{0in}This work was completed while D. Dominici was visiting Technische
Universit\"{a}t Berlin and supported in part by a Sofja Kovalevskaja Award
from the Humboldt Foundation, provided by Professor Olga Holtz. He wishes to
thank Olga for her generous sponsorship and his colleagues at TU Berlin for
their continuous help.
\end{acknowledgement}

\end{document}